\documentclass [12pt]{amsart}
\usepackage{amssymb,amsmath,amsthm}

\newtheorem{lemma}{Lemma}
\newtheorem{proposition}{Proposition}

\newcommand{\tr}{\,\mathrm{tr}\,}
\newcommand{\ad}{\,\mathrm{ad}\,}

\newcommand{\GL}{\,\mathrm{GL}\,}
\newcommand{\SL}{\,\mathrm{SL}\,}
\newcommand{\Aut}{\,\mathrm{Aut}\,}

\newcommand{\diag}{\,\mathrm{diag}\,}
\begin{document}

\begin{center}

{\Large {\bf Automorphisms of elementary adjoint Chevalley groups\\

\bigskip

of types $A_l, D_l, E_l$ over local rings with $1/2$\footnote{The
work is supported by the Russian President grant MK-2530.2008.1 and
by the grant of Russian Fond of Basic Research 08-01-00693.} }}

\bigskip
\bigskip

{\large \bf E.~I.~Bunina}

\end{center}
\bigskip

\begin{center}

{\bf Abstract.}

\end{center}

In this paper we prove that every automorphism of and elementary
adjoint Chevalley group of types $A_l, D_l$, or $E_l$, over local
commutative ring with $1/2$ is a composition of a ring automorphism
and conjugation by some matrix from the normalizer of the Chevalley
group in $\GL(V)$ ($V$ is the adjoint representation space).

\bigskip

\section*{Introduction}\leavevmode


Let $G_{\ad}$ be a Chevalley--Demazure group scheme associated with
an irreducible root system~$\Phi$ of type $A_l$ ($l\ge 2$), $D_l$
($l\ge 4$), $E_l$ ($l=6,7,8$); $G_{\ad}(\Phi,R)$ be a set of
points~$G_{\ad}$ with values in a commutative ring~$R$;
$E_{\ad}(\Phi,R)$ be the elementary subgroup of~$G_{\ad}(\Phi,R)$,
where $R$ is a commutative ring with~$1$. In this paper we describe
automorphisms of groups $E_{\ad}(\Phi,R)$ over local commutative
rings with~$1/2$.

Similar results for Chevalley groups over fields were proved by
R.\,Steinberg~\cite{Stb1} for finite case and by
J.\,Humphreys~\cite{H} for infinite case. Many papers were devoted
to description of automorphisms of Chevalley groups over different
commutative rings, we can mention here the papers of
Borel--Tits~\cite{v22}, Carter--Chen~Yu~\cite{v24},
Chen~Yu~\cite{v25}--\cite{v29}. E.\,Abe~\cite{Abe_OSN} proved that
automorphisms are standard for Noetherian rings, it could completely
close the question about automorphisms of Chevalley groups over
arbitrary commutative rings (for the case of systems of rank $\ge 2$
and rings with $1/2$), but in consideration of adjoint elementary
groups in the paper~\cite{Abe_OSN} there is a mistake, that can not
be corrected by methods of this paper. Namely, in the proof of
Lemma~11 the author uses the fact that $\mathrm{ad}\,
(x_\alpha)^2=0$ for all long roots, but it is not true in the
adjoint representation. The main problem here is the case of groups
of type~$E_8$, since in all other cases Chevalley groups have a
representation with the property $\mathrm{ad}\,(x_\alpha)^2=0$ for
all long roots, but in the case~$E_8$ there are no such
representations.

In the given paper we consider also this  case $E_8$, and it can
help to close the question about automorphisms of Chevalley groups
over commutative rings with~$1/2$.

 We generalize some methods of V.M.~Petechuk~\cite{Petechuk1} to prove the main theorem.

Note that we consider the cases $A_l$, $D_l$, $E_l$, but the case
$A_l$ was completely studied by the papers of
W.C.~Waterhouse~\cite{v46}, V.M.~Petechuk~\cite{v12},  Fuan Li and
Zunxian Li~\cite{v37}, and also for rings without~$1/2$. The paper
of I.Z.\,Golubchik and A.V.~Mikhalev~\cite{v8} covers the
case~$C_l$, that is not considered in the present paper.

The author is thankful to N.A.\,Vavilov, A.Yu.\,Golubkov,
A.A.\,Klyachko, A.V.\,Mikhalev, and specially to the referee of this
paper for valuable advices, remarks and discussions.

\section{Definitions and formulation of the main theorem}\leavevmode

We fix a root system~$\Phi$, that has one of types $A_l$ ($l\ge 2$),
$D_l$ ($l\ge 4$), or $E_l$ ($l=6,7,8$), with the system of simple
roots~$\Delta$, the set of positive (negative) roots $\Phi^+$
($\Phi^-$), the Weil group~$W$. Recall that in our case any two
roots are conjugate under the action of the Weil group. Let
$|\Phi^+|=m$. More detailed texts about root systems and their
properties can be found in the books \cite{Hamfris}, \cite{Burbaki}.

Suppose now that we have a semisimple complex Lie algebra~$\mathcal
L$ with the Cartan subalgebra~$\mathcal H$ (more details about
semisimple Lie algebras can be found in the book~\cite{Hamfris}).

Lie algebra  $\mathcal L$ has a decomposition ${\mathcal
L}={\mathcal H} \oplus \sum\limits_{\alpha\ne 0} {\mathcal
L}_\alpha$,
$$
{\mathcal L}_\alpha:=\{ x\in {\mathcal L}\mid [h,x]=\alpha(h)x\text{
for every } h\in {\mathcal H}\},
$$
and if ${\mathcal L}_\alpha\ne 0$, then $\dim {\mathcal
L}_\alpha=1$, all nonzero $\alpha\in {\mathcal H}$ such that
${\mathcal L}_\alpha\ne 0$, form some root system~$\Phi$. The root
system $\Phi$ and the semisimple Lie algebra
 $\mathcal L$ over~$\mathbb C$
uniquely (up to automorphism) define each other.

On the Lie algebra $\mathcal L$ we can introduce a bilinear
\emph{Killing form} $\varkappa(x,y)=\tr (\ad x\ad y),$  that is
non-degenerated on~$\mathcal H$. Therefore we can identify the
spaces $\mathcal H$ and ${\mathcal H}^*$.

We can choose a basis $\{ h_1, \dots, h_l\}$ in~$\mathcal H$ and for
every $\alpha\in \Phi$ elements $x_\alpha \in {\mathcal L}_\alpha$
so that $\{ h_i; x_\alpha\}$ is a basis in~$\mathcal L$ and for
every two elements of this basis their commutator is an integral
linear combination of the elements of the same basis. This basis is
called a \emph{Chevalley basis}.

Introduce now elementary Chevalley groups (see~\cite{Steinberg}).

Let  $\mathcal L$ be a semisimple Lie algebra (over~$\mathbb C$)
with a root system~$\Phi$, $\pi: {\mathcal L}\to \mathfrak{gl}(V)$
be its finitely dimensional faithful representation  (of
dimension~$n$). If $\mathcal H$ is a Cartan subalgebra of~$\mathcal
L$, then a functional
 $\lambda \in {\mathcal H}^*$ is called a
 \emph{weight} of  a given representation, if there exists a nonzero vector $v\in V$
 (that is called a  \emph{weight vector}) such that
for any $h\in {\mathcal H}$ $\pi(h) v=\lambda (h)v.$

In the space~$V$ in the Chevalley basis all operators
$\pi(x_\alpha)^k/k!$ for $k\in \mathbb N$ are written as integral
(nilpotent) matrices. An integral matrix also can be considered as a
matrix over an arbitrary commutative ring with~$1$. Let $R$ be such
a ring. Consider matrices $n\times n$ over~$R$, matrices
$\pi(x_\alpha)^k/k!$ for
 $\alpha\in \Phi$, $k\in \mathbb N$ are included in $M_n(R)$.

Now consider automorphisms of the free module $R^n$ of the form
$$
\exp (tx_\alpha)=x_\alpha(t)=1+t\pi(x_\alpha)+t^2
\pi(x_\alpha)^2/2+\dots+ t^k \pi(x_\alpha)^k/k!+\dots
$$
Since all matrices $\pi(x_\alpha)$ are nilpotent, we have that this
series is finite. Automorphisms $x_\alpha(t)$ are called
\emph{elementary root elements}. The subgroup in $\Aut(R^n)$,
generated by all $x_\alpha(t)$, $\alpha\in \Phi$, $t\in R$, is
called an \emph{elementary Chevalley group} (notation:
$E_\pi(\Phi,R)$).

In elementary Chevalley group we can introduce the following
important elements and subgroups:

--- $w_\alpha(t)=x_\alpha(t) x_{-\alpha}(-t^{-1})x_\alpha(t)$, $\alpha\in \Phi$,
$t\in R^*$;

--- $h_\alpha (t) = w_\alpha(t) w_\alpha(1)^{-1}$;

---  $N$ is generated by all
 $w_\alpha (t)$, $\alpha \in \Phi$, $t\in R^*$;

---  $H$ is generated by all
 $h_\alpha(t)$, $\alpha \in \Phi$, $t\in R^*$;

The action of  $x_\alpha(t)$ on the Chevalley basis is described in
\cite{v23}, \cite{VavPlotk1}.

It is known that the group $N$ is a normalizer of~$H$ in elementary
Chevalley group, the quotient group $N/H$ is isomorphic to the Weil
group $W(\Phi)$.

All weights of a given representation (by addition) generate a
lattice (free Abelian group, where every  $\mathbb Z$-basis  is also
a $\mathbb C$-basis in~${\mathcal H}^*$), that is called the
\emph{weight lattice} $\Lambda_\pi$.

 Elementary Chevalley groups are defined not even by a representation of the Chevalley groups,
but just by its \emph{weight lattice}. Namely, up to an abstract
isomorphism an elementary Chevalley group is completely defined by a
root system~$\Phi$, a commutative ring~$R$ with~$1$ and a weight
lattice~$\Lambda_\pi$.

Among all lattices we can mark two: the lattice corresponding to the
adjoint representation, it is generated by all roots (the \emph{root
lattice}~$\Lambda_{ad}$) and the lattice generated by all weights of
all reperesentations (the \emph{lattice of weights}~$\Lambda_{sc}$).
For every faithful reperesentation~$\pi$ we have the inclusion
$\Lambda_{ad}\subseteq \Lambda_\pi \subseteq \Lambda_{sc}.$
Respectively, we have the \emph{adjoint} and \emph{universal}
elementary Chevalley groups. In this paper we study adjoint
elementary Chevalley groups.

Every elementary Chevalley group satisfies the following conditions:

(R1) $\forall \alpha\in \Phi$ $\forall t,u\in R$\quad
$x_\alpha(t)x_\alpha(u)= x_\alpha(t+u)$;

(R2) $\forall \alpha,\beta\in \Phi$ $\forall t,u\in R$\quad
 $\alpha+\beta\ne 0\Rightarrow$
$$
[x_\alpha(t),x_\beta(u)]=x_\alpha(t)x_\beta(u)x_\alpha(-t)x_\beta(-u)=
\prod x_{i\alpha+j\beta} (c_{ij}t^iu^j),
$$
where $i,j$ are integers, product is taken by all roots
$i\alpha+j\beta$, replacing in some fixed order; $c_{ij}$ are
integer numbers not depending of $t$ and~$u$, but depending of
$\alpha$ and $\beta$ and of order of roots in the product. In the
cases under consideration always
$$
[x_\alpha(t),x_\beta(u)]=x_{\alpha+\beta}(\pm tu).
$$

(R3) $\forall \alpha \in \Phi$ $w_\alpha=w_\alpha(1)$;

(R4) $\forall \alpha,\beta \in \Phi$ $\forall t\in R^*$ $w_\alpha
h_\beta(t)w_\alpha^{-1}=h_{w_\alpha (\beta)}(t)$;

(R5) $\forall \alpha,\beta\in \Phi$ $\forall t\in R^*$ $w_\alpha
x_\beta(t)w_\alpha^{-1}=x_{w_\alpha(\beta)} (ct)$, where
$c=c(\alpha,\beta)= \pm 1$;

(R6) $\forall \alpha,\beta\in \Phi$ $\forall t\in R^*$ $\forall u\in
R$ $h_\alpha (t)x_\beta(u)h_\alpha(t)^{-1}=x_\beta(t^{\langle
\beta,\alpha \rangle} u)$.

By $X_\alpha$ we denote the subgroup consisting of all $x_\alpha
(t)$ for $t\in R$.

We will need two types of automorphisms of elementary Chevalley
groups
 $E_{\ad}(\Phi,R)$.

{\bf Ring automorphisms.} Let $\rho: R\to R$ be an automorphism
of~$R$. The mapping $x\mapsto \rho (x)$ from $E_{\ad}(\Phi,R)$ onto
itself is an automorphism of the group $E_{\ad}(\Phi,R)$, tha is
denoted by the same letter~$\rho$ and is called a \emph{ring
automorphism} of the group~$G_\pi(\Phi,R)$. Note that for all
$\alpha\in \Phi$ and $t\in R$ an element $x_\alpha(t)$ is mapped
into $x_\alpha(\rho(t))$.

{\bf Automorphisms-conjugations.} Let $V$ be a reprepresentation
space of the group $E_{\ad} (\Phi,R)$, $C\in \GL(V)$ be some matrix
that does not move our Chevalley group:
$$
C E_{\ad}(\Phi,R) C^{-1}= E_{\ad} (Phi,R).
$$
 Then the mapping $x\mapsto CxC^{-1}$ from $E_\pi(\Phi,R)$ onto itself is an automorphism of the Chevalley group,
 that is denoted by $i_С$ and is called
an \emph{automorphism-conjugation} of the group~$E(R)$,
\emph{induced by the element}~$C$ of~$\GL(V)$.

\bigskip

{\bf  Theorem 1.} \emph{Let $E_{\ad}(\Phi,R)$ be an elementary
Chevalley group with an irreducible root system of type $A_l$ $(l\ge
2)$, $D_l$ $(l\ge 4)$, or $E_l$ $(l=6,7,8)$, $R$ be a commutative
local ring with~$1/2$. Then every automorphism of $E_{\ad}(\Phi,R)$
is a composition of a ring automorphism and an
automorphism-conjugation.}

Next sections are devoted to the proof of Theorem~1.

\section{Replacing the initial automorphism to the special one.}\leavevmode

From this section we suppose that $R$ is a local ring with $1/2$,
the Chevalley group is adjoint, in this section root system is
arbitrary. In this section we use some reasonings
from~\cite{Petechuk1}.

Let $J$ be the maximal ideal (radical) of~$R$, $k$  the residue
field $R/J$. Then $E_J=E_{ad}( \Phi,R,J)$ is the greatest normal
proper subgroup of $E_{\ad}(\Phi,R)$ (see~\cite{Abe1}). Therefore,
$E_J$ is invariant under the action of~$\varphi$.

By this reason  the automorphism
$$
\varphi: E_{\ad} (\Phi,R)\to E_{\ad}(\Phi,R)
$$
induces an automorphism
$$
\overline \varphi: E_{\ad} (\Phi,R)/E_J=E_{\ad} (\Phi,k)\to
E_{\ad}(\Phi,k).
$$
The group $E_{\ad}(\Phi,k)$ is a Chevalley group over field,
therefore the automorphism $\overline \varphi$ is standard (see
\cite{Steinberg}), i.\,e. it has the form
$$
\overline \varphi =    i_{\overline g} \overline \rho,\quad
\overline g\in N(E_{\ad}(\Phi,k)),
$$
where $\overline \rho$ is a ring automorphism, induced by some
automorphism of~$k$.

 It is clear that there exists a matrix $g\in GL_n(R)$ such that
its image under factorization  $R$ by~$J$ coincides with~$\overline
g$. We are not sure that $g\in N(E_{\ad}(\Phi,R))$.

Consider  a mapping $\varphi'= i_{g^{-1}} \varphi$. It is an
isomorphism of the group
 $E_{ad}(\Phi,R)\subset GL_n(R)$ onto some subgroup in $GL_n(R)$,
with the property that its image under factorization $R$ by $J$
coincides with the automorphism $\overline \rho$.

These arguments prove

\begin{proposition}\label{dop1}
Every matrix $A\in E_{\ad}(\Phi,R)$ with elements from the
subring~$R'$ of~$R$, generated by unit, is mapped under the action
of~$\varphi'$ to some matrix from the set
$$
A\cdot \GL_n(R,J)=\{ B\in \GL_n(R)\mid A-B\in M_n(J)\}.
$$
\end{proposition}

Let $a\in E_{ad} (\Phi,R)$, $a^2=1$. Then the element $e=\frac{1}{2}
(1+a)$ is an idempotent in the ring $M_n(R)$. This idempotent  $e$
defines a decomposition of the free $R$-module $V=R^n$:
$$
V=eV\oplus (1-e)V=V_0\oplus V_1
$$
(the modules $V_0$, $V_1$ are free, since every projective module
over local field is free~\cite{Mc}). Let $\overline V=\overline V_0
\oplus \overline V_1$ be decomposition of the $k$-module~$\overline
V$ with respect to~$\overline a$, and
 $\overline e=\frac{1}{2} (1+\overline a)$.

Then we have

\begin{proposition}\label{pr1_1}
The modules \emph{(}subspaces\emph{)}
 $\overline V_0$, $\overline V_1$ are images of the modules $V_0$, $V_1$ under factorization by~$J$.
\end{proposition}
\begin{proof} Let us denote the images of $V_0$, $V_1$ under factorization
by $J$ by $\widetilde V_0$, $\widetilde V_1$, respectively. Since
$V_0=\{ x\in V| ex=x\},$ $V_1= \{ x\in V|ex=0\},$  we have
 $\overline e(\overline x)=\frac{1}{2}(1+\overline a)(\overline x)=\frac{1}{2}
(1+\overline a(\overline
x))=\frac{1}{2}(1+\overline{a(x)})=\overline{e(x)}$. Then
$\widetilde V_0\subseteq \overline V_0$, $\widetilde V_1\subseteq
\overline V_1$.

Let $x=x_0+x_1$, $x_0\in V_0$, $x_1\in V_1$. Then $\overline
e(\overline x)=\overline e(\overline x_0)+\overline e (\overline
x_1)=\overline x_0$. If $\overline x\in \widetilde V_0$, then
$\overline x=\overline x_0$.
\end{proof}

 Let
$b=\varphi'(a)$. Then $b^2=1$ and $b$ is equivalent to $a$
modulo~$J$.

\begin{proposition}\label{pr1_2}
Suppose that $ a,b\in E_\pi(\Phi,R)$, $a^2=b^2=1$, $a$ is a matrix
with elements from the subring of~$R$, generated by the unit, $b$
and $a$ are equivalent modulo~$J$, $V=V_0\oplus V_1$ is a
decomposition of~$V$ with respect to~$a$, $V=V_0'\oplus V_1'$ is a
decomposition of~$V$ with respect to~$b$. Then $\dim V_0'=\dim V_0$,
$\dim V_1'=\dim V_1$.
\end{proposition}

\begin{proof}
We have an $R$-basis of the module~$V$ $\{ e_1,\dots,e_n\}$ such
that $\{ e_1,\dots,e_k\}\subset V_0$, $\{ e_{k+1},\dots,
e_n\}\subset V_1$.  It is clear that
$$
\overline a \overline e_i=\overline{ae_i}=\overline {(\sum_{j=1}^n
a_{ij} e_j)}= \sum_{j=1}^n \overline a_{ij} \overline e_j.
$$
Let $\overline V=\overline V_0\oplus \overline V_1$, $\overline V =
\overline V_0'\oplus \overline V_1'$ are decompositions of
$k$-module (space)~$\overline V$ with respect to $\overline a$ and $
\overline b$. It is clear that $\overline V_0= \overline V_0'$,
$\overline V_1 =\overline V_1'$. Therefore, by
Proposition~\ref{pr1_1} the images of the modules $V_0$ and $V_0'$,
$V_1$ and $V_1'$ under factorization by~$J$ coincide. Let us take
such $\{ f_1,\dots, f_k\}\subset V_0'$, $\{ f_{k+1},\dots,
f_n\}\subset V_1'$ that $\overline f_i=\overline e_i$,
$i=1,\dots,n$. Since the matrix of transformation from $\{
e_1,\dots, e_n\}$ to $\{ f_1,\dots, f_n\}$ is invertible (it is
equivalent to the identical matrix modulo~$J$) we have that $\{
f_1,\dots, f_n\}$ is a $R$-basis in~$V$. It is clear that $\{
f_1,\dots, f_k\}$ is a $R$-basis in $V_0'$, $\{ v_{k+1},\dots,
v_n\}$ is a $R$-basis in $V_1'$.
\end{proof}

\section{Images of~$w_{\alpha_i}$}

We consider some fixed adjoint Chevalley group $E=E_{ad}(\Phi,R)$
with the root system $A_l$ ($l\ge 2$), $D_l$ ($l\ge 4$), $E_6$,
$E_7$ or $E_8$, its adjoint representation in the group $GL_n(R)$
($n=l+2m$, where $m$ is the number of positive roots of~$\Phi$),
with the basis of weight vectors $v_1=x_{\alpha_1},
v_{-1}=x_{-\alpha_1}, \dots, v_n=x_{\alpha_n}, v_{-n}=x_{-\alpha_n},
V_1=h_{1},\dots,V_l=h_{l}$, corresponding to the Chevalley basis of
the system~$\Phi$.

We also have the isomorphism~$\varphi'$, described in Section~2.

Consider the matrices $h_{\alpha_1}(-1),\dots, h_{\alpha_l}(-1)$ in
our basis. They have the form
$$
h_{\alpha_i}(-1)=\diag [\pm 1,\dots, \pm 1,
\underbrace{1,\dots,1}_{l}],
$$
on $(2j-1)$-th and $(2j)$-th places we have $-1$ if and only if
 $\langle \alpha_i,\alpha_j\rangle=-1$. As we see, for all~$i$
$h_{\alpha_i}(-1)^2=1$.

According to Proposition~\ref{pr1_2} we know that every matrix
$h_i=\varphi''(h_{\alpha_i}(-1))$ in some basis is diagonal with
$\pm 1$ on the diagonal, and the number of $1$ and $-1$ coincides
with its number for the matrix $h_{\alpha_i}(-1)$. Since all
matrices $h_i$ commutes, there exists a basis, where all $h_i$ have
the same form as $h_{\alpha_i}(-1)$. Suppose that we come to this
basis with the help of the matrix~$g_1$. It is clear that
 $g_1\in GL_n(R,J)$. Consider the mapping
 $\varphi_1=i_{g_1}^{-1} \varphi'$. It is also an isomorphism
of the group $E$ onto some subgroup of $GL_n(R)$ such that its image
under factorization $R$ by~$J$ is~$\overline \rho$, and
$\varphi_1(h_{\alpha_i}(-1))=h_{\alpha_i}(-1)$ for all
$i=1,\dots,l$.

Let us consider the isomorphism~$\varphi_1$.

Matrices
$h_{\alpha_k}(a)=diag[a_1,1/a_1,a_2,1/a_2,\dots,a_m,1/a_m,1,\dots,1]$
commute with all $h_{\alpha_i}(-1)$, therefore their images under
the isomorphism~$\varphi_1$ also commute with
all~$h_{\alpha_i}(-1)$. Thus, these images have the form
$$
\begin{pmatrix}
C_1& 0& \dots &0& 0\\
0& C_2& \dots& 0& 0\\
\vdots& \vdots& \ddots& \vdots& \vdots\\
0& 0& \dots& C_n& 0\\
0& 0& \dots& 0& C
\end{pmatrix},\quad
C_i\equiv \begin{pmatrix}
\overline \rho(a_i)& 0\\
0& \overline \rho(1/a_i)
\end{pmatrix} \ \mod J,\quad C\in \GL_l(R, J).
$$

Every element $w_i=w_{\alpha_i}(1)$ maps (by conjugation) diagonal
matrices into diagonal ones, so its image has block-monomial form.

From $\varphi_1(w_i)\equiv w_i\mod J$ it follows that blocks of
$\varphi_1(w_i)$ are in the same places as blocks of~$w_i$.

Consider the first  vector of the basis obtained after the last
change. Let us denote it by~$e$. The Weil group $W$ transitively
acts on all  roots, therefore for every root~$\alpha_i$ there exists
such $w^{(\alpha_i)}\in W$ that $w^{(\alpha_i)} \alpha_1=\alpha_i$.
 Consider now
the basis $e_1,\dots, e_{2m}, e_{2m+1},\dots, e_{2m+l}$, where
$e_1=e$, $e_i=\varphi_1(w^{(\alpha_i)})e$; for $2m< i\le 2m+1$ $e_i$
is not changed. It is clear that the matrix of this basis change is
equivalent to $1$ modulo~$J$. Therefore the obtained set of vectors
is a basis.

It is clear that the matrix $\varphi_1(w_i)$ ($i=1,\dots,l$) in the
part of basis $\{ e_1,\dots,e_{2m}\}$
 coincides with the matrix $w_i$ in the initial basis of weight vectors.
Since $h_i(-1)$ are squares of $w_i$, then their images also are not
changed in the new basis.

Besides that, we know that  $\varphi_1(w_i)$ is block-diagonal up to
the first $2m$ and last $l$ elements. Therefore, the last basis
part, consisting of $l$ elements, can be changed independently.

Let us denote matrices $w_i$ and $\varphi_1(w_i)$ on this part of
basis by $\widetilde w_i$ and $\widetilde{\varphi_1(w_i)}$,
respectively. All these matrices are involutions, they have only one
$-1$ in their diagonal forms. Let $\widetilde V=\widetilde
V_0^i\oplus \widetilde V_1^i$ be decomposition of the matrix
$\widetilde{\varphi_1(w_i)}$.

\begin{lemma}\label{l3_1}
Matrices $\widetilde{\varphi_1(w_i)}$ and
$\widetilde{\varphi_1(w_j)}$, where $i\ne j$, commute if and only if
$\widetilde V_1^i\subseteq \widetilde V_0^j$ and $\widetilde
V_1^j\subseteq \widetilde V_0^i$.
\end{lemma}
\begin{proof}
If $\widetilde{\varphi_1(w_i)}$ and $\widetilde{\varphi_1(w_j)}$
commute, then the (free one-dimensional) submodule $\widetilde
V_1^i$ is proper for $\widetilde{\varphi_1(w_j)}$ and the (free
one-dimensional) submodule $\widetilde V_1^j$ is proper for
$\widetilde{\varphi_1(w_i)}$. Therefore either $\widetilde
V_1^i\subset \widetilde V_1^j$ or $\widetilde V_1^i\subset
\widetilde V_0^j$. If $\widetilde V_1^i\subset \widetilde V_1^j$
then $\widetilde V_1^i=\widetilde V_1^j$. Since the module $V_0^i$
is invariant for $\widetilde \varphi_1(w_j)$, we have $\widetilde
V_0^i\subset \widetilde V_0^j$, therefore $\widetilde
V_0^i=\widetilde V_0^j$, and so $\widetilde
\varphi_1(w_i)=\widetilde \varphi_1(w_j)$ and we come to
contradiction. Consequently, $\widetilde V_1^i\subset \widetilde
V_0^j$, and similarly $\widetilde V_1^j\subset \widetilde
V_0^i$.\end{proof}

\begin{lemma}\label{l3_2}
For any root system $\Phi$ there exists such a basis in $\widetilde
V$ that the matrix $\widetilde \varphi_1(w_1)$ in this basis has the
same form as $w_1$, i.e. is equal to
$$
\begin{pmatrix}
-1& 1& 0\\
0& 1& 0\\
0& 0& E_{l-2}
\end{pmatrix}.
$$
\end{lemma}
\begin{proof}
Since $\widetilde w_1$ is an involution and $\widetilde V_1^1$ has
dimension $1$, there exists a basis $\{ e_1,e_2,\dots, e_l\}$ where
$\widetilde \varphi_1(w_1)$ has the form $diag [-1,1,\dots,1]$. In
the basis $\{ e_1, e_2-1/2e_1,e_3,\dots, e_l\}$ the matrix
$\widetilde \varphi_1(w_1)$ has the obtained form.
\end{proof}

\begin{lemma}\label{l3_3}
For the root system $A_2$ there exists such a basis that $\widetilde
\varphi_1(w_1)$ and $\widetilde \varphi_1(w_2)$ in this basis have
the same form as $w_1$ and $w_2$, i.e. are equal to
$$
\begin{pmatrix}
-1& 1\\
0& 1
\end{pmatrix} \text{ and }
\begin{pmatrix}
1& 0\\
1& -1
\end{pmatrix},
$$
respectively.
\end{lemma}
\begin{proof}
By Lemma~\ref{l3_2} we can find a basis in~$\widetilde V$ such that
 the matrix $\widetilde \varphi_1(w_1)$ in this basis has the same
form as $w_1$. Let the matrix $\widetilde \varphi_1(w_2)$ in this
basis be
$$
\begin{pmatrix}
a& b\\
c& d
\end{pmatrix}.
$$
Let us make the basis change with the help of the matrix
$$
\begin{pmatrix}
c& (1-c)/2\\
0& 1
\end{pmatrix}.
$$
Under this basis change the matrix $\widetilde \varphi_1(w_1)$
remains the same form, and the matrix $\widetilde \varphi_1(w_2)$
becomes
$$
\begin{pmatrix}
a'& b'\\
1& d'
\end{pmatrix}.
$$
As this matrix is an involution, we have $a'+d'=0$, ${a'}^2+b'=1$.
So we obtain $d'=-a'$. Now let us use the condition
$$
\left(\begin{pmatrix} -1& 1\\
0& 1
\end{pmatrix} \begin{pmatrix}
a'& b'\\
1& -a'
\end{pmatrix}\right)^2 =
\begin{pmatrix}
a'& b'\\
1& -a'
\end{pmatrix}
\begin{pmatrix} -1& 1\\
0& 1
\end{pmatrix}.
$$
This condition gives (its second line and first row) $1-2a'=-1$,
therefore $a'=1$. From ${a'}^2+b'=1$ it follows $b'=0$.
\end{proof}

\begin{lemma}\label{l3_4}
For every root system $\Phi\ne A_2$ we can choose a basis in
$\widetilde V$ such that the matrices $\widetilde \varphi_1(w_1)$
and  $\widetilde \varphi_1(w_2)$ in this basis have the same form as
$\widetilde w_1$ and $\widetilde w_2$, respectively.
\end{lemma}
\begin{proof}
The intersection of modules $\widetilde V_0^1$ and $\widetilde
V_0^2$ is a free module of dimension $\ge l-3$. Therefore we can
suppose that $\widetilde \varphi_1(w_1)$ and  $\widetilde
\varphi_1(w_2)$ have the  form $\begin{pmatrix} *& 0\\
0& E_{l-3} \end{pmatrix}$. Moreover, by Lemma~\ref{l3_2} we can
suppose that $\widetilde \varphi_1(w_1)$  has the same form as
$\widetilde w_1$. We can consider not the whole module~$\widetilde
V$, but its limitation to the first three basis vectors. Let
$$
\widetilde \varphi_1(w_1)=\begin{pmatrix} a_1& a_2& a_3\\
b_1& b_2& b_3\\
c_1& c_2& c_3
\end{pmatrix}.
$$
Taking basis change with the matrix
$$
\begin{pmatrix}
b_1& (1-b_1)/2& 0\\
0& 1& 0\\
0& 0& 1
\end{pmatrix},
$$
we do not change $\widetilde \varphi_1(w_1)$, but $\widetilde
\varphi_1(w_2)$ becomes
$$
\widetilde \varphi_1(w_1)=\begin{pmatrix} a_1'& a_2'& a_3'\\
1& b_2'& b_3'\\
c_1'& c_2'& c_3'
\end{pmatrix}.
$$
Now we will use the same conditions as in the previous lemma. The
first is $\widetilde \varphi_1(w_2)^2-1=0$ (Cond. 1) and the second
is $(\widetilde \varphi_1(w_1)\widetilde
\varphi_1(w_2))^2-\widetilde \varphi_1(w_2)\widetilde
\varphi_1(w_1)=0$ (Cond. 2). If we subtract Condition 1 from
Condition 2 we obtain (line 2, row 1) $a_1'=1$, (line 2, row 2)
$a_2'=0$, then from Cond. 1, line 1, row 3, we obtain
$a_3'(1+c_3')=0$. As $c_3'\equiv 1\mod J$, we have $a_3'=0$. The
same condition, line 2, row 3, gives $b_3'(b_2'+c_3')=0$, as
$b_3'\in R^*$, we have $c_3'=-b_2'$.

Again taking basis change, but  with the matrix
$$
\begin{pmatrix}
1& 0& 0\\
0& 1& 0\\
0& -c_1'& 1
\end{pmatrix},
$$
we do not change $\widetilde \varphi_1(w_1)$, but $\widetilde
\varphi_1(w_2)$ becomes
$$
\widetilde \varphi_1(w_1)=\begin{pmatrix} 1& 0& 0\\
1& b_2''& b_3''\\
0& c_2''& -b_2''
\end{pmatrix}.
$$
Then directly from Cond.~1 we obtain $b_2''=-1$, $c_2''=0$, and the
last basis change with the matrix $diag[1,1,b_3'']$ makes the
obtained forms of  $\widetilde \varphi_1(w_1)$  and $\widetilde
\varphi_1(w_2)$.
\end{proof}

\begin{lemma}\label{l3_5}
For the root system $D_4$ there exists such a basis that $\widetilde
\varphi_1(w_1)$, $\widetilde \varphi_1(w_2)$, $\widetilde
\varphi_1(w_3)$ and $\widetilde \varphi_1(w_4)$ in this basis have
the same forms as $\widetilde w_1$, $\widetilde w_2$, $\widetilde
w_3$, and $\widetilde w_4$, i.e. are equal to
$$
\begin{pmatrix}
-1& 1& 0& 0\\
0& 1& 0& 0\\
0& 0& 1& 0\\
0& 0& 0&1
\end{pmatrix}, \begin{pmatrix}
1& 0& 0& 0\\
1& -1& 1& 1\\
0& 0& 1& 0\\
0& 0& 0& 1 \end{pmatrix}, \begin{pmatrix} 1& 0& 0& 0\\
0& 1& 0& 0\\
0& 1& -1& 1\\
0& 0& 0& 1 \end{pmatrix} \text{ and }
\begin{pmatrix}
1& 0& 0& 0\\
0& 1& 0& 0\\
0& 0& 1& 0\\
0& 0& 1& -1
\end{pmatrix},
$$
respectively.
\end{lemma}
\begin{proof}
We take such a basis that $\widetilde \varphi_1(w_1)$, $\widetilde
\varphi_1(w_3)$, $\widetilde \varphi_1(w_4)$ have the same form as
the initial $\widetilde w_1$, $\widetilde w_3$, $\widetilde w_4$. We
can do it because $w_1$, $w_3$, $w_4$ are commuting involutions,
there exists a basis where  $\widetilde \varphi_1(w_1)$, $\widetilde
\varphi_1(w_3)$, $\widetilde \varphi_1(w_4)$ have the forms
$diag[-1,1,1,1]$, $diag[1,1,-1,1]$, $diag[1,1,1,-1]$, respectively.
Then, conjugating them by the matrix $$ \begin{pmatrix} 1& -1/2& 0&
0\\
0& 1& 0& 0\\
0& -1/2& 1& 0\\
0& -1/2& 0& 1 \end{pmatrix},
 $$
 we come to the obtained basis. Now let us look for $\widetilde \varphi_1(w_2)$. We have the
following conditions: $\widetilde w_2^2=1$ (Cond.~1), $(\widetilde
w_1 \widetilde w_2)^2=\widetilde w_2 \widetilde w_1$ (Cond.~2),
$(\widetilde w_3\widetilde w_2)^2=\widetilde w_2\widetilde w_3$
(Cond.~3), $(\widetilde w_4\widetilde w_2)^2=\widetilde
w_2\widetilde w_4$ (Cond.~4). Let
$$
\widetilde \varphi_1(w_2)=\begin{pmatrix} a_1& a_2& a_3& a_4\\
b_1& b_2& b_3& b_4\\
c_1& c_2& c_3& c_4\\
d_1& d_2& d_3& d_4
\end{pmatrix}.
$$

Taking basis change with the matrix
$$
\begin{pmatrix}
1& 0& 0& 0\\
0& 1& 0& 0\\
0& 0& 1& 0\\
0& \frac{d_1}{2d_1-b_1} & 0& \frac{b_1}{b_1-2d_1}
\end{pmatrix},
$$
we do not change $\widetilde \varphi_1(w_1)$, $\widetilde
\varphi_1(w_3)$, $\widetilde \varphi_1(w_4)$, but $\widetilde
\varphi_1(w_2)$ becomes
$$
\begin{pmatrix}
a_1& a_2& a_3& a_4\\
b_1& b_2& b_3& b_4\\
c_1& c_2& c_3& c_4\\
0& d_2& d_3& d_4
\end{pmatrix}
$$
(we do not write primes for simplicity).

Now from the 4-th line of (Cond.~1$-$ Cond.~2) it follows
$d_2=d_3=0$, $d_4=1$.

Line 2, row 4 of (Condition 1 $-$ Condition 4) gives $b_4(b_4-1)=0$.
Since $b_4\in R^*$, we have $b_4=1$.

Now, taking basis change with the matrix
$$
\begin{pmatrix}
\frac{b_3}{b_3-2a_3}& \frac{a_3}{2a_3-b_3}& 0& 0\\
0& 1& 0& 0\\
0& 0& 1& 0\\
0& 0 & 0& 1
\end{pmatrix},
$$
we do not change $\widetilde \varphi_1(w_1)$, $\widetilde
\varphi_1(w_3)$, $\widetilde \varphi_1(w_4)$, but $\widetilde
\varphi_1(w_2)$ becomes
$$
\begin{pmatrix}
a_1& a_2& 0& a_4\\
b_1& b_2& b_3& 1\\\
c_1& c_2& c_3& c_4\\
0& 0& 0& 1
\end{pmatrix}
$$
(again we do not write primes for simplicity).

Then line 1, row 3 of Cond.~1 gives $a_2b_3=0\Rightarrow a_2=0$,
line 1, row 1 of Cond.~1 gives $a_1^2=1\Rightarrow a_1=1$, line 1,
row 4 gives $2a_4=0\Rightarrow a_4=0$. Line 2, row 4 of
Cond.~1$-$Cond.~2 gives $b_1=1$, line 3, row~4 gives $c_4=c_1$. Line
2, row 3 of Cond.1 gives $c_3=-b_2$.

Finally, taking basis change with the matrix
$$
\begin{pmatrix}
1& 0& 0& 0\\
0& 1& 0& 0\\
0& \frac{1-b_3}{2}& b_3& 0\\
0& 0 & 0& 1
\end{pmatrix},
$$
we do not change $\widetilde \varphi_1(w_1)$, $\widetilde
\varphi_1(w_3)$, $\widetilde \varphi_1(w_4)$, but $\widetilde
\varphi_1(w_2)$ becomes
$$
\begin{pmatrix}
1& 0& 0& a_4'\\
1& b_2'& 1& 1\\\
c_1'& c_2'& -b_2'& c_1'\\
0& 0& 0& 1
\end{pmatrix}.
$$
After that line 2, row 3 of Cond.~3 gives $b_2'=-1$, and then
Cond.~1 gives $c_1'=c_2'=0$.
\end{proof}

\begin{lemma}\label{l3_6}
Suppose that we have some root system $\Phi$ and elements
$\widetilde \varphi_1(w_{i_1})=\widetilde w_{i_1}$,\dots,
$\widetilde \varphi_1(w_{i_k})=\widetilde w_{i_k}$, and also
$\widetilde \varphi_1(w_{i_{k+1}})$, where one of the following
cases holds:

a) $\Phi=A_l$, $l\ge 3$, $i_1=1$, $i_2=2$,\dots, $i_k=k$,
$i_{k+1}=k+1$, $k+1< l$;

b) $\Phi = D_l$, $l> 4$, $i_1=l$, $i_2=l-1$, $i_3=l-2$, \dots ,
$i_k=l-k+1$, $i_{k+1}=l-k$, $4< k< l$;

c) $\Phi=E_6, E_7$ or $E_8$, $i_1=1$, $i_2=2$,\dots, $i_k=k$,
$i_{k+1}=k+1$, $4\le k< l-1$.

Then we can choose such a basis of $\widetilde V$ that $\widetilde
\varphi_1(w_{i_1})=\widetilde w_{i_1}$,\dots, $\widetilde
\varphi_1(w_{i_{k+1}})=\widetilde w_{i_{k+1}}$.
\end{lemma}
\begin{proof}
In all these cases $\widetilde \varphi_1(w_{i_{k+1}})$ commutes with
all $\widetilde \varphi_1(w_{i_1})=\widetilde w_{i_1}$,\dots,
$\widetilde \varphi_1(w_{i_{k-1}})=\widetilde w_{i_{k-1}}$,
therefore (see Lemma~\ref{l3_1}) for all $j=i_1,\dots ,i_{k-1}$
$V_1^j\subset V_0^{i_{k+1}}$. Since $V_1^{i_1}\oplus \dots \oplus
V_1^{i_{k-1}}=\langle e_{i_1},\dots , e_{i_{k-1}}\rangle$, we infer
that $\widetilde \varphi_1(w_{i_{k+1}})$ is identical on the first
$k-1$ basic vectors. As in Lemma~\ref{l3_4} we obtain that
$\widetilde \varphi_1(w_{i_{k+1}})$ is identical on the last $l-k-2$
basic vectors. Therefore we can limit $\widetilde
\varphi_1(w_{i_{k+1}})$ for the part of basis $\{ e_k, e_{k+1},
e_{k+2}\}$ (without loss of generality). Now the proof is completely
the same as in Lemma~\ref{l3_4}.
\end{proof}

\begin{proposition}\label{pr3_1}
For every root system $\Phi=A_l, D_l, E_l$ we can choose a basis in
$\widetilde V$ such that the matrices $\widetilde
\varphi_1(w_1)$,\dots, $\widetilde \varphi_1(w_l)$ in this basis
have the same form as $\widetilde w_1$,\dots, $\widetilde w_l$,
respectively.
\end{proposition}
\begin{proof}
If we have the system $A_2$, we can use Lemma~\ref{l3_3}. If
$\Phi=A_l$, $l\ge 3$, then we apply Lemma~\ref{l3_4}, after that
Lemma~\ref{l3_6} $l-3$ times, and finally the same arguments as in
Lemma~\ref{l3_3} for the element   $\widetilde \varphi_1(w_l)$.

If $\Phi=D_l$, then we apply Lemma~\ref{l3_5} to the roots
$\alpha_{l-3},\dots, \alpha_l$, then Lemma~\ref{l3_6} $l-5$ times to
the roots $\alpha_{l-4},\dots, \alpha_2$, and finally the same
arguments as in Lemma~\ref{l3_3} for the element   $\widetilde
\varphi_1(w_1)$.

If $\Phi=E_l$, then we apply Lemma~\ref{l3_5} to the roots
$\alpha_{2},\dots, \alpha_5$, then Lemma~\ref{l3_6}   to the roots
$\alpha_{6},\dots, \alpha_{l-1}$, and finally the same arguments as
in Lemma~\ref{l3_3} for the elements   $\widetilde \varphi_1(w_1)$
and $\widetilde \varphi_1(w_l)$.
\end{proof}

Therefore, we can now consider the isomorphism $\varphi_2$ with all
properties of $\varphi_1$, and such that $\varphi_2(w_i)=w_i$ for
all $i=1,\dots,l$.

We suppose that we have the isomorphism~$\varphi_2$ with these
properties.

\section{The images of $x_{\alpha_i}(1)$ and diagonal matrices.}

Now we are interested in the images of $x_{\alpha_i}(t)$. Let
$\varphi_2(x_{\alpha_1}(1))=x_1$. Since $x_1$ commutes with all
$h_{\alpha_i}(-1)$, $i=1,3,\dots, l$, we have that  $x_1$ is
separated to the blocks of the following form: blocks $2\times 2$
are corresponded to the part of basis $\{ v_i,v_{-i}\}$, where $i>
1$ and
 $\langle \alpha_i,\alpha_1\rangle \ne 0$;
blocks $4\times 4$ are corresponded to the part of basis
 $\{ v_i,v_{-i},v_j,v_{-j}\}$,
where $i> 1$, $\alpha_i=\alpha_j\pm \alpha_1$; and we also have the
part
 $\{ v_1,v_{-1},V_1,\dots, V_l\}$.

For $h_{\alpha_2}(-1)$ we know that $h_{\alpha_2}(-1) x_1
h_{\alpha_2}(-1)= x_1^{-1}$.  Then, on the blocks $2\times 2$,
described above if $x_1$ has the form
$$
\begin{pmatrix}
a& b\\
c& d
\end{pmatrix},
$$
then this matrix to the second power is $1$, therefore
$a^2+bc=d^2+bc=1$, $b(a+d)= c(a+d)=0$. Since $a+d\equiv 2\mod J$, we
have $a+d\in R^*$, i.\,e. $b=c=0$. Since $a^2=d^2=1$ and $a,d\equiv
1\mod J$, we have $a=d=1$. Thus, on the blocks $2\times 2$ the
matrix $x_1$ is always $1$ (i.e. it coincides with $x_{\alpha_1}(t)$
on these blocks), so we can now not to consider these basis
elements.

Now let us use the conditions $w_i x_1 w_i^{-1} =x_1$ for $i\ge 3$.

At first, all blocks $4\times 4$ has the same form, since every two
such blocks are conjugate up to the action of $w_i$, $i\ge 3$.

The conditions $w_i x_1 w_i^{-1} =x_1$ for $i\ge 3$ for the rest of
the basis together with the condition $h_2 x_1 h_2^{-1} =x_1^{-1}$
say that the matrix on the basis subset $\{ v_1,v_{-1},
V_1,V_2,V_3,\dots,V_l\}$  has the form
$$
\begin{pmatrix}
* & * & * & *& 0&\dots&  0\\
* & * & * & *& 0&\dots& 0\\
* & * & * & *& 0&\dots& 0\\
* & * & * & *& 0&\dots& 0\\
* & * & * & *&1& \dots& 0\\
\hdotsfor{7}\\
*& * & * & * & 0& \dots& 1
\end{pmatrix},
$$
and all lines $5,\dots,l$ are expressed via the fourth line.
According to the zero corner of this matrix we can restrict the
 conditions to its left upper submatrix $4\times 4$.

Suppose that on this part of the basis the matrix $x_1$ has the form
$$
\begin{pmatrix}
a_1& a_2& a_3& a_4\\
b_1& b_2& b_3& b_4\\
c_1& c_2& c_3& c_4\\
d_1& d_2& d_3& d_4
\end{pmatrix},
$$
and on the part of basis $v_2,v_{-2},v_{1+2},v_{-1-2}$ it has the
form
 $$ \begin{pmatrix} e_1& e_2& e_3& e_4\\
f_1& f_2& f_3& f_4\\
g_1& g_2& g_3& g_4 \\
h_1& h_2& h_3& h_4
\end{pmatrix}. $$

We will consider the part of basis $\{ v_1,v_{-1},v_2,v_{-2},
v_{1+2}, v_{-1-2}, V_1,V_2\}$.

Taking basis change with the block-diagonal matrix, that has the
form
$$
\begin{pmatrix}
1& -\frac{b_4}{a_4}\\
-\frac{b_4}{a_4}& 1
\end{pmatrix}
$$
on every block $\{ v_i, v_{-i}\}$ (it is possible, because $b_4\in
J$), and is identical on the block
 $\{ V_1,\dots, V_l\}$, we do not change
the elements  $w_i, h_i$, and $x_1$ now has $b_4=0$. So we can
suppose that the isomorphism $\varphi_2$ is such that
$\varphi_2(x_{\alpha_1}(1))$ has $b_4=0$.

Then we make the basis change with the help  of diagonal matrix,
having the form $\frac{1}{a_4}\cdot E$ on the part  $\{
v_1,v_{-1},\dots, v_m,v_{-m}\}$, and being identical on the part $\{
V_1,\dots, V_l\}$. Similarly, all elements $w_i,h_i$ are not
changed, and $a_4$ is now equal to~$1$.

So we suppose that $\varphi_2(x_{\alpha_1}(1))$ has $b_4=0$ and
$a_4=1$.

On the part of basis under consideration
$$
w_1=\begin{pmatrix}
0& -1& 0& 0& 0& 0& 0& 0\\
-1& 0& 0& 0& 0& 0& 0& 0\\
0& 0& 0& 0& 1& 0& 0& 0\\
0& 0& 0& 0& 0& 1& 0& 0\\
0& 0& -1& 0& 0& 0& 0& 0\\
0& 0& 0& -1& 0& 0& 0& 0\\
0& 0& 0& 0& 0& 0& -1& 1\\
0& 0& 0& 0& 0& 0& 0& 1
\end{pmatrix},\quad
w_2=\begin{pmatrix}
0& 0& 0& 0& -1& 0& 0& 0\\
0& 0& 0& 0& 0& -1& 0& 0\\
0& 0& 0& -1& 0& 0& 0& 0\\
0& 0& -1& 0& 0& 0& 0& 0\\
1& 0& 0& 0& 0& 0& 0& 0\\
0& 1& 0& 0& 0& 0& 0& 0\\
0& 0& 0& 0& 0& 0& 1& 0\\
0& 0& 0& 0& 0& 0& 1& -1
\end{pmatrix}.
$$
Since
$$
x_{\alpha_1}(1)=\begin{pmatrix}
1& -1& 0& 0& 0& 0& -2& 1\\
0& 1& 0& 0& 0& 0& 0& 0\\
0& 0& 1& 0& 0& 0& 0& 0\\
0& 0& 0& 1& 0& 1& 0& 0\\
0& 0& -1& 0& 1& 0& 0& 0\\
0& 0& 0& 0& 0& 1& 0& 0\\
0& 1& 0& 0& 0& 0&1& 0\\
0& 0& 0& 0& 0& 0& 0& 1
\end{pmatrix},
$$
we have
$$
x_1=\varphi_2(x_1(1))=\begin{pmatrix} a_1&a_2&0& 0& 0& 0& a_3& 1\\
b_1& b_2& 0& 0& 0& 0& b_3& 0\\
0& 0& e_1& e_2& e_3& e_4& 0& 0\\
0& 0& f_1& f_2& f_3& f_4& 0& 0\\
0& 0& g_1& g_2& g_3& g_4& 0& 0\\
0& 0& h_1& h_2& h_3& h_4& 0& 0\\
c_1& c_2& 0& 0&0& 0& c_3& c_4\\
d_1& d_2& 0& 0& 0&0& d_3& d_4
\end{pmatrix},
$$
where $a_1,b_2,e_1, f_2,f_4, g_3, h_4, c_2, c_3, d_4\equiv 1\mod J$,
$a_2,g_1\equiv -1\mod J$, $a_3\equiv -2\mod J$, all other entries
are in~$J$.

 Then
$$
x_{1+2}=\varphi_2(x_{\alpha_1+\alpha_2}(1))=w_2x_1w_2^{-1}=\begin{pmatrix}
g_3&
g_4& g_2& g_1& 0& 0& 0& 0\\
h_3& h_4& h_2& h_1& 0& 0& 0& 0\\
f_3& f_4& f_2& f_1&0& 0&0&0\\
e_3&e_4& e_2& e_1& 0& 0& 0& 0\\
0&0& 0& 0& a_1& a_2& 1+a_3& -1\\
0& 0& 0& 0& b_1& b_2& b_3& 0\\
0& 0& 0& 0& c_1& c_2& c_3+c_4& -c_4\\
0& 0& 0& 0& c_1-d_1& c_2-d_2& c_3-d_3+c_4-d_4& -c_4+d_4
\end{pmatrix}
$$
and
$$
x_2=\varphi_2(x_{\alpha_2}(1))=w_1x_{1+2}w_1^{-1}=\begin{pmatrix}
h_4& h_3&
0& 0& h_2& h_1& 0& 0\\
g_4& g_3& 0& 0&g_2& g_1& 0& 0\\
0& 0& a_1& a_2& 0& 0& -1-a_3& a_3\\
0& 0& b_1& b_2& 0& 0& -b_3& b_3\\
f_4& f_3& 0& 0& f_2& f_1& 0& 0\\
e_4& e_3& 0& 0& e_2& e_1& 0& 0\\
0& 0& -d_1& -d_2& 0& 0& d_3+d_4& -d_3\\
0& 0& c_1-d_1& c_2-d_2& 0& 0& -c_3+d_3-c_4+d_4& c_3-d_3
\end{pmatrix}.
$$
We will use the following conditions:
$$
\mathrm{Con1}:=(x_1x_{1+2}-x_{1+2}x_1=0),\quad \mathrm{Con2}:=(h_2
x_1 h_2 x_1-1=0).
$$
Position (3,8) of $\mathrm{Con1}$ gives $f_3=-e_3$, position (2,8)
of $\mathrm{Con1}$ gives $h_3=-b_3c_4$, position (2,8) of
$\mathrm{Con2}$ gives $b_1=b_3c_4$, therefore $h_3=-b_1$. From
position (1,1) of $\mathrm{Con1}$ we have $b_1(a_2+g_4)=0$,
therefore $b_1=0$, since $a_2+g_4\in R^*$.

Now we introduce two more conditions:
$$
\mathrm{Con3}:=(x_1w_1x_1w1^{-1}-w_1h_2x_1h_2=0),\quad
\mathrm{Con4}:=(x_2x_1-x_{1+2}x_1x_2=0).
$$

Denote $y_1=a_1-1, y_2=a_2+1, y_3=a_3+2, y_4=b_2-1, y_5=b_3,
y_6=c_1, y_7=c_2-1, y_8=c_3-1, y_9=c_4, y_{10}=d_1, y_{11}=d_2,
y_{12}=d_3, y_{13}=d_4-1, y_{14}=e_1-1, y_{15}=e_2, y_{16}=e_3$,
$y_{17}=e_4, y_{18}=f_1, y_{19}=f_2-1, y_{20}=f_4-1, y_{21}=g_1+1,
y_{22}=g_2, y_{23}=g_3-1, y_{24}=g_4, y_{25}=h_1, y_{26}=h_2,
y_{27}=h_4-1$. All these $y_i$ are in~$J$. From the Conditions 1--4
we have the following $27$ equations (that are linear with respect
to $y_i$):

\begin{align*}
&y_{23}(-a_2)+y_{24}(a_1-b_2)+y_{27}a_2=0,& \text{ Con1,
pos.~(1,2)},\\
& y_{18}(-g_1)+y_{22}(a_1-e_1)+y_{26}(a_2)=0,&\text{ Con1,
pos.~(1,3)},\\
&y_1g_1+y_{15}(-g_2)+y_{19}(-g_1)+y_{25}a_2=0,& \text{ Con1,
pos.~(1,4)},\\
&y_6(a_3+1)+y_{10}(-1)+y_{16}(g_1+g_2)=0,& \text{ Con1,
pos.~(1,5)},\\
&y_3(c_2)+y_7(-1)+y_{11}(-1)+y_{20}+y_{21}(-f_4)+y_{22}(-e_4)=0,&
\text{ Con1, pos.~(1,6)},\\
&y_3(c_3+c_4-g_3)+y_8(-1)+y_9(-1)+y_{13}(-1)+y_{14}(-1)+y_{23}2+&\\
&+y_{24}(-b_3)=0,&\text{
Con1, pos.~(1,7)},\\
&y_9(-a_3-1)+y_{13}+y_{23}(-1)=0,&\text{ Con1, pos.~(1,8)},\\
&y_5c_2+y_{25}(-f_4)+y_{26}(-e_4)=0,&\text{ Con1, pos.~(2,6)},\\
&y_{16}a_2+y_{17}(b_2-f_2)+y_{18}(-f_4)=0,& \text{ Con1,
pos.~(3,6)},\\
&y_5(e_4-f_4)+y_{16}(1+2a_3)=0,&\text{ Con1, pos.~(3,7)},\\
&y_{15}(-f_1f_4-f_2e_3)+y_{16}(a_1-a_2+a_1h_2+f_1f_2-f_2^2)+y_{22}(e_3a_2-f_4b_2)=0,&\text{
Con4, pos.~(3,5)},\\
&y_{10}(-d_3-d_4)+y_{11}(a_1+1)+y_{12}c_1=0,& \text{ Con3,
pos.~(8,2)},\\
&y_1(-1)+y_2(a_1+b_2)+y_3(-c_2)+y_4(-1)+y_72+y_{11}(-1)=0,&\text{
Con2, pos.~(1,2)},\\
&y_5(-c_2g_3)+y_{16}(b_2-b_2h_4)+y_{17}a_2=0,&\text{ Con4,
pos.~(6,2)},\\
&y_{14}g_3+y_{16}(e_4-e_3)+y_{21}+y_{23}=0,& \text{ Con3,
pos.~(3,3)},\\
&y_4(b_2+1)+y_5(-c_2)=0,&\text{ Con2, pos.~(2,2)},\\
&y_{14}(e_1+1)+y_{15}f_1+y_{16}(-g_1)+y_{17}(-h_1)=0,& \text{ Con2,
pos.~(3,3)},\\
&y_{15}(e_1+f_2)+y_{16}(-g_2)+y_{17}(-h_2)=0,&\text{ Con2,
pos.~(3,4)},\\
&y_6(-g_1a_1)+y_9(c_3+c_4-d_4)(c_1-d_1)+y_{10}(c_3^2+c_3c_4-d_3c_4-e_1)+&\\
&+y_{11}(-f_1)+y_{25}c_2a_1=0,&\text{
Con4, pos.~(7,3)},\\
&y_{16}g_4+y_{18}e_4+y_{19}+y_{20}(f_2-h_4)+y_{27}(-1)=0,&\text{
Con2, pos.~(4,6)},\\
&y_{14}+y_{21}(g_3-e_1)+y_{22}f_1+y_{23}(-1)+y_{24}h_1=0,&\text{
Con2, pos.~(5,3)},\\
&y_{17}(-g_1)+y_{22}(-f_4)+y_{24}(g_3+h_4)=0,& \text{ Con2,
pos.~(5,6)},\\
&y_4(-c_2)+y_6(-a_2)+y_8c_2+y_9d_2=0,&\text{ Con2, pos.~(7,2)},\\
&y_6(-1)+y_9(c_3+d_4)=0,&\text{ Con2, pos.~(7,8)},\\
&y_1a_2+y_2+y_4+y_6a_3+y_{10}(-1-a_3)=0,&\text{ Con3, pos.~(1,2)},\\
&y_{19}h_4+y_{20}(-1)+y_{25}(-g_2)+y_{26}(-h_2)+y_{27}=0,&\text{
Con3, pos.~(6,6)},\\
&y_6(-g_3f_2)+y_{15}(-c_1g_4-c_2h_4)+y_{16}(d_2-d_1)+&\\
&+y_{22}(c_4d_2-c_3c_2-c_2c_4)+y_{26}(c_4d_1-c_3c_1-c_4c_1)=0,&\text{
Con4, pos.~(7,5)}.
\end{align*}

The matrix of this system of linear equations modulo~$J$ is  {\tiny
$$
 \left(\begin{array}{ccccccccccccccccccccccccccc} 0& 0& 0& 0& 0&
0& 0& 0& 0& 0& 0& 0& 0& 0& 0& 0& 0&0 &0& 0& 0& 0& 1&
0&0 & 0& -1\\
0& 0& 0& 0& 0& 0& 0& 0& 0& 0& 0& 0& 0& 0& 0& 0& 0&0 &1& 0& 0& 0& 0&
0&0 & -1& 0\\
-1& 0& 0& 0& 0& 0& 0& 0& 0& 0& 0& 0& 0& 0& 0& 0& 0&0 &1& 0& 0& 0& 0&
0&-1 & 0& 0\\
0& 0& 0& 0& 0& -1& 0& 0& 0& -1& 0& 0& 0& 0& 0& -1& 0&0 &0& 0& 0& 0&
0& 0&0 & 0& 0\\
0& 0& 1& 0& 0& 0& -1& 0& 0& 0& -1& 0& 0& 0& 0& 0& 0&0 &0& 1& -1& 0&
0& 0&0 & 0& 0\\
0& 0& 0& 0& 0& 0& 0& 1& 1& 0& 0& 1& 1& 0& 0& 0& 0&0 &0& 0& 0& 0&
-2& 0&0 & 0& 0\\
0& 0& 0& 0& 0& 0& 0& 0& 1& 0& 0& 0& 1& 0& 0& 0& 0&0 &0& 0& 0& 0& -1&
0&0 & 0& 0\\
0& 0& 0& 0& 1& 0& 0& 0& 0& 0& 0& 0& 0& 0& 0& 0& 0&0 &0& 0& 0& 0& 0&
0&-1 & 0& 0\\
0& 0& 0& 0& 0& 0& 0& 0& 0& 0& 0& 0& 0& 0& 0& 1& 0&1 &0& 0& 0& 0&
0& 0&0 & 0& 0\\
0& 0& 0& 0& -1& 0& 0& 0& 0& 0& 0& 0& 0& 0& 0& -3& 0&0 &0& 0& 0& 0&
0& 0&0 & 0& 0\\
0& 0& 0& 0& 0& 0& 0& 0& 0& 0& 0& 0& 0& 0& 0& 1& 0&0 &0& 0& 0& -1& 0&
0&0 & 0& 0\\
0& 0& 0& 0& 0& 0& 0& 0& 0& -1& 2& 0& 0& 0& 0& 0& 0&0 &0& 0& 0& 0& 0&
0&0 & 0& 0\\
-1& 2& -1& -1& 0& 0& 2& 0& 0& 0& -1& 0& 0& 0& 0& 0& 0&0 &0& 0& 0& 0&
0&0&0 & 0& 0\\
0& 0& 0& 0& 1& 0& 0& 0& 0& 0& 0& 0& 0& 0& 0& 0& 1&0 &0& 0& 0& 0&
0&0&0 & 0& 0\\
0& 0& 0& 0& 0& 0& 0& 0& 0& 0& 0& 0& 0& 1& 0& 0& 0&0 &0& 0& 1& 0& 1&
0&0 & 0& 0\\
0& 0& 0& 2& -1& 0& 0& 0& 0& 0& 0& 0& 0& 0& 0& 0& 0&0 &0& 0& 0& 0& 0&
0&0 & 0& 0\\
0& 0& 0& 0& 0& 0& 0& 0& 0& 0& 0& 0& 0& 2& 0& 1& 0&0 &0& 0& 0& 0& 0&
0&0 & 0& 0\\
0& 0& 0& 0& 0& 0& 0& 0& 0& 0& 0& 0& 0& 0& 2& 0& 0&0 &0& 0& 0& 0& 0&
0&0 & 0& 0\\
0& 0& 0& 0& 0& 1& 0& 0& 0& 0& 0& 0& 0& 0& 0& 0& 0&0 &0& 0& 0& 0& 0&
0&1 & 0& 0\\
0& 0& 0& 0& 0& 0& 0& 0& 0& 0& 0& 0& 0& 0& 0& 0& 0&0 &1& 0& 0& 0& 0&
0&0 & 0& -1\\
0& 0& 0& 0& 0& 0& 0& 0& 0& 0& 0& 0& 0& 1& 0& 0& 0&0 &0& 0& 0& 0& -1&
0&0 & 0& 0\\
0& 0& 0& 0& 0& 0& 0& 0& 0& 0& 0& 0& 0& 0& 0& 0& 1&0 &0& 0& 0& -1& 0&
2&0 & 0& 0\\
0& 0& 0& -1& 0& 1& 0& 1& 0& 0& 0& 0& 0& 0& 0& 0& 0&0 &0& 0& 0& 0& 0&
0&0 & 0& 0\\
0& 0& 0& 0& 0& -1& 0& 0& 2& 0& 0& 0& 0& 0& 0& 0& 0&0 &0& 0& 0& 0& 0&
0&0 & 0& 0\\
-1& 1& 0& 1& 0& -2& 0& 0& 0& 1& 0& 0& 0& 0& 0& 0& 0&0 &0& 0& 0& 0&
0&0&0 & 0& 0\\
0& 0& 0& 0& 0& 0& 0& 0& 0& 0& 0& 0& 0& 0& 0& 0& 0&0 &1& -1& 0& 0& 0&
0&0 & 0& 1\\
0& 0& 0& 0& 0& 1& 0& 0& 0& 0& 0& 0& 0& 0& 1& 0& 0&0 &0& 0& 0& 1&
0&0&0 & 0& 0
\end{array}\right).
$$
}

Determinant of this matrix is $2^8$, so it is invertible in~$R$.
Therefore this system has the unique solution $y_1=\dots=y_{27}=0$.
Consequently, $x_1=x_{\alpha_1}(1)$ on the part of basis under
consideration. Since all roots are conjugate up to the action
of~$W$, we have that $x_1=x_{\alpha_1}(1)$ on the whole basis. It is
clear that also $x_2=x_{\alpha_2}(1)$.

  Now consider the
matrix $d_t=\varphi_1(h_{\alpha_1}(t))$. The matrix
$h_{\alpha_1}(t)$ is $diag[t^2,1/t^2,1/t,t,t,1/t,1,1]$ on the part
of basis under consideration.

\begin{lemma}\label{l4_1}
The matrix $d_t$ is $h_{\alpha_1}(s)$ for some $s\in R^*$.
\end{lemma}

\begin{proof}
For the matrix $d_t$ we have conditions  $d_t w_3=w_3 d_t$, \dots,
$d_t w_l =w_l d_t$.

Let $l> 2$ and
$$
d_t=\begin{pmatrix} \gamma_{11} & \gamma_{12} & \dots &
\gamma_{1l}\\
\gamma_{21}& \gamma_{22}& \dots& \gamma_{2l}\\
\vdots& \vdots& \ddots& \vdots\\
\gamma_{l1}& \gamma_{l2}& \dots & \gamma_{ll}
\end{pmatrix}
$$
on~$\widetilde V$.

Каждое из соотношений $d_t w_i=w_i d_t$, $i>2$, дает
$\gamma_{1i}=\dots=\gamma_{i-1,i}=\gamma_{i+1,i}=\dots=\gamma_{li}=0$.
Из соотношения $d_t w_1 d_t w_1^{-1}=1$ теперь сразу следует
$\gamma_{33}^2=\dots=\gamma_{ll}^2=1\Rightarrow
\gamma_{33}=\dots=\gamma_{ll}=1$. Условие $d_t w_l=w_l d_t$
показывает, что $\gamma_{l-1,j}$ линейно выражается через
$\gamma_{l,j}$, $j=1,2$, \dots, условие $d_t w_3=w_3 d_t$
показывает, что $\gamma_{2,j}$ линейно выражается через
$\gamma_{l,j}$, $j=1,2$.

Every condition $d_t w_i=w_i d_t$, $i>2$, gives
$\gamma_{1i}=\dots=\gamma_{i-1,i}=\gamma_{i+1,i}=\dots=\gamma_{li}=0$.
From the condition $d_t w_1 d_t w_1^{-1}=1$ now directly follows
$\gamma_{33}^2=\dots=\gamma_{ll}^2=1\Rightarrow
\gamma_{33}=\dots=\gamma_{ll}=1$. Condition $d_t w_l=w_l d_t$ gives
that $\gamma_{l-1,j}$ is linearly expressed via $\gamma_{l,j}$,
$j=1,2$,  \dots, condition $d_t w_3=w_3 d_t$ gives that
$\gamma_{2,j}$ is linearly expressed via $\gamma_{l,j}$, $j=1,2$.

According to the zero corner we can consider the conditions for
$d_t$ on the part of basis $v_1,v_{-1},v_2,v_{-2},v_{1+2}$,
$v_{-1-2},V_1,V_2$, as we did it for $x_1$.

Since $d_t$ commutes with all $h_i$, we have that on the part of
basis under consideration
$$
d_t=\begin{pmatrix} k_1&k_2& 0& 0& 0& 0& 0& 0\\
k_3& k_4& 0& 0& 0& 0& 0& 0\\
0& 0& l_1& l_2& 0& 0& 0& 0\\
0& 0& l_3& l_4& 0& 0& 0& 0\\
0& 0& 0& 0& m_1& m_2& 0& 0\\
0& 0& 0& 0& m_3& m_4& 0& 0\\
0& 0& 0& 0& 0& 0& n_1& n_2\\
0& 0& 0& 0& 0& 0& n_3& n_4
\end{pmatrix}.
$$

We know that $x_2=x_{\alpha_2}(1)$ and
$x_2^t=\varphi_2(x_{\alpha_2}(t))=d_tx_2d_t^{-1}=d_tx_2w_1d_tw_1^{-1}$.
Using these conditions we obtain
the expression of $x_2^t$ trough the entries of $d_t$. Then we can
use the condition $\mathrm{Con5}:=(x_2^t x_2-x_2x_2^t=0)$. Its
position (1,6) gives $k_1(k_2+k_3)=0\Rightarrow k_2=-k_3$,
pos.~(2,1) gives $k_4(l_3-m_3)=0\Rightarrow m_3=l_3$, pos.~(2,5)
gives $l_3(l_1+m_4)=0\Rightarrow l_3=0$, pos.~(5,6) gives
$k_2(l_4+m_1)=0\Rightarrow k_2=0$. From
$\mathrm{Con6}:=(d_tw_1d_tw_1^{-1}-1=0)$ it follows
$k_1k_4=l_1m_1=l_4m_4=1$. Using the condition
$\mathrm{Con7}:=(w_2d_tw_2^{-1}-d_tw_1w_2d_tw_2^{-1}w_1^{-1}=0)$, we
obtain $m_2=l_2=0$ (positions (1,2) and (4,3)) and $l_4=l_1k_1$
(pos.~(1,1)). Position (7,7) of $\mathrm{Con6}$ gives
$n_1^2-(n_1+n_2)n_3=1$, position (3,7) of $\mathrm{Con5}$ gives
$n_1^2-(3n_1+n_2-2n_3-n_4)n_3=1$. Therefore, $n_3=0$, $n_1^2=1$.
After that we clearly infer $n_1=n_4=1$, $n_2=0$. Finally, position
(3,4) of $\mathrm{Con5}$ gives $k_1=1/l_1^2$. Denote $1/l_1$ by~$s$.

It is clear that with the help of the elements
 $w_i$, $i=3,\dots l$, we can define all other diagonal elements.
Namely, if $\langle \alpha_1, \alpha_k\rangle=p$, then
$\varphi(h_{\alpha_1}(t)) v_k = s^p\cdot v_k$,
$\varphi(h_{\alpha_1}(t)) v_{-k} = s^{-p}\cdot v_{-k}$. Whence
$\varphi(h_{\alpha_1}(t))=h_{\alpha_1}(s)$.
\end{proof}

\section{Images of  the matrices $x_{\alpha_i}(t)$,
proof of the main theorem.}

It is clear that $\varphi_2(h_{\alpha_k}(t))=h_{\alpha_k}(s)$,
$k=1,\dots,n$. Let us denote the mapping $t\mapsto s$ by $\rho: R^*
\to R^*$. Note that for $t\in R^*$
$\varphi_2(x_1(t))=\varphi_2(h_{\alpha_2}(t^{-1}) x_1(1)
h_{\alpha_2}(t))=h_{\alpha_2}(s^{-1}) x_1(1) h_{\alpha_2}(s)=
x_1(s)$. If $t\notin R^*$, then $t\in J$, i.\,e. $t=1+t_1$, where
$t_1\in R^*$. Then
$\varphi_2(x_1(t))=\varphi_2(x_1(1)x_1(t_1))=x_1(1)x_1(\rho(t_1))=
x_1(1+\rho(t_1))$. Therefore, if we extend the mapping $\rho$ on the
whole ring~$R$ (by the formula $\rho(t):=1+\rho(t-1)$, $t\in R$), we
have $\varphi_2(x_1(t))=x_1(\rho(t))$ for all $t\in R$. It is clear
that $\rho$ is injective, additive, multiplicative on all invertible
elements. Since every element of~$R$ is the sum of two invertible
elements, we have that  $\rho$ is an isomorphism from~$R$ onto some
its subring~$R'$. Note that in this situation $C E(\Phi,R)
C^{-1}=E(\Phi,R')$ for some matrix $C\in GL(V)$. Let us show that
$R'=R$.

Let us denote the matrix units by $E_{ij}$.

\begin{lemma}\label{porozhd}
Elementary Chevalley group $E(\Phi,R)$ generates the ring $M_n(R)$.
\end{lemma}
\begin{proof}
The matrix $(x_{\alpha_1}(1)-1)^2$ has the single nonzero element
$-2\cdot E_{12}$. Multiplying it to some suitable diagonal matrix we
can obtain an arbitrary matrix of the form $\lambda\cdot E_{12}$
(since $-2\in R^*$ and $R^*$ generates~$R$). According to the
transitive action of the Weil group on the root system~$\Phi$ (for
every root $\alpha_k$ there exists such an element $w\in W$ that
$w(\alpha_1)=\alpha_k$) the matrix $\lambda E_{12}\cdot w$ has the
form $\lambda E_{1,2k}$,  and the matrix $w^{-1}\cdot \lambda
E_{12}$ has the form $\lambda E_{2k-1,2}$. Moreover, according to
the Weil group element that maps the first root to the opposite one,
we get an element $E_{2,1}$. Taking different combinations of the
obtained elements, we can get an arbitrary element $\lambda E_{ij}$,
$1\le i,j\le 2m$. Therefore we have always generated the subring
 $M_{2m}(R)$.
Now let us subtract from $x_{\alpha_1}(1)-1$ suitable matrix units,
and we obtain the matrix $E_{2m+1, 2}-2 E_{1,2m+1}+ E_{1,2m+2}$.
Multiplying its (from the right side) to $E_{2,i}$, $1\le i\le 2m$,
we get all $E_{2m+1, i}$, $1\le i\le 2m$. With the help of Weil
group elements we have all $E_{i,j}$, $2m< i\le 2m+l$, $1\le j\le
2n$. Now we have the matrix $-2E_{1,2m+1}+E_{1,2m+2}$. Multiplying
it (from the left side) to $E_{2m+1,1}$, we get $E_{2m+1,2m+1}$.
With the help of two last matrices we have $E_{1,2m+1}$, and,
therefore, all
 $E_{i,j}$, $1\le i\le 2m$, $2m< j\le 2m+l$.
It is clear that now we have all matrix units, i.e. the whole matrix
ring $M_n(R)$.

Let us show it for the simplest root system $A_2$. In this case
$(x_{\alpha_1}(1)-1)^2=-2E_{12}$, $h_{\alpha_2}(t) (-2E_{12})=-2t
E_{12}$, i.e., we can obtain every $\lambda E_{12}$. Then
$w_{\alpha_1} \lambda E_{12} w_{\alpha_1}(1)^{-1} =\lambda E_{21}$,
$\lambda E_{12}E_{21}=\lambda E_{11}$, $\lambda E_{21}E_{12}=\lambda
E_{22}$, $w_{\alpha_2}(1) \lambda E_{12}=\lambda E_{52}$,
$w_{\alpha_2}(1) \lambda E_{21}=\lambda E_{61}$, $\lambda E_{52}
E_{21}=\lambda E_{51}$, $\lambda E_{61} E_{12}=\lambda E_{62}$,
$\lambda E_{12} w_{\alpha_2}(1)=\lambda E_{16}$, $\lambda E_{21}
w_{\alpha_2}(1)=\lambda E_{25}$, $\lambda E_{21} E_{16} =\lambda
E_{26}$, $\lambda E_{12} E_{25}=\lambda E_{15}$, $\lambda E_{51}
E_{15}=\lambda E_{55}$, $\lambda E_{61}E_{16}=\lambda E_{66}$,
$\lambda E_{51}E_{16}=E_{56}$, $\lambda E_{61}E_{15}=E_{65}$,
$\lambda E_{i5} w_{\alpha_1}(1)=\lambda E_{i3}$, $i=1,2,5,6$,
$\lambda E_{i6} w_{\alpha_1}(1)=\lambda E_{i4}$, $i=1,2,5,6$,
$\lambda w_{\alpha_1}(1) E_{5i}=\lambda E_{3i}$, $i=1,2,5,6$,
$\lambda w_{\alpha_1}(1) E_{6i}=\lambda E_{4i}$, $i=1,2,5,6$,
$\lambda E_{41} E_{13}=\lambda E_{43}$, $\lambda E_{41}
E_{14}=\lambda E_{44}$, $\lambda E_{31} E_{13}=\lambda E_{33}$,
$\lambda E_{31}E_{14}=\lambda E_{34}$, so we have all matrix units
of the subring $M_6(R)$.

Then
$y=x_{\alpha_1}(1)-1=-E_{12}-2E_{17}+E_{18}+E_{46}-E_{53}+E_{73}$,
$y'=y+E_{12}-E_{46}+E_{53}=E_{18}-2E_{17}+E_{72}$,
$(E_{18}-2E_{17}+E_{72})\cdot \lambda E_{2i}=\lambda E_{7i}$,
$i=1,\dots,6$, $(w_{\alpha_2}(1)-1) \lambda E_{7i} = \lambda
E_{8i}$, $i=1,\dots,6$, $y''=y'-E_{72}=E_{18}-2E_{17}$, $\lambda
E_{81} y''=\lambda E_{88}$, $\lambda E_{71} y'' =-2\lambda E_{77}$,
$y'' \lambda E_{88} =\lambda E_{18}$, $y'' \lambda E_{77}=-2 \lambda
E_{17}$, $\lambda E_{i1} E_{17} =\lambda E_{i7}$, $\lambda E_{i1}
E_{18} =\lambda E_{i8}$, so we have generated the whole ring
$M_8(R)$.
\end{proof}

\begin{lemma}\label{Tema}
If for some $C\in \GL(V)$ we have $C E(\Phi,R) C^{-1}= E(\Phi,R')$,
where $R'$ is a subring of~$R$, then $R'=R$.
\end{lemma}
\begin{proof}
Suppose that $R'$ is a proper subring of~$R$.

Then $C M_n(R) C^{-1} =M_n (R')$, since the group $E(\Phi,R)$
generates the ring $M_n(R)$, and the group $E(\Phi,R')=CE(\Phi,R)
C^{-1}$ generates the ring $M_n(R')$. It is impossible, since $C\in
GL_n(R)$.
\end{proof}

Now we have proved that  $\rho$ is an automorphism of the ring~$R$.
Therefore, composition of the initial automorphism~$\varphi$ and
some basis change with the help of the matrix $C\in \GL_n(R)$ (that
maps $E(\Phi,R)$ onto itself), is a ring automorphism~$\rho$. It
proves Theorem 1. $\square$

\end{document}